\documentclass[reqno]{amsart}
\usepackage{amssymb}
\usepackage[colorlinks=true,hypertex]{hyperref}
\date{This manuscript was finalized on 16 June 2006}
\date{}
\allowdisplaybreaks[4]
\theoremstyle{plain}
\newtheorem{thm}{Theorem}
\DeclareMathOperator{\td}{d\mspace{-1mu}}
\newcommand{\abs}[1]{\bigl\lvert#1\bigr\rvert}

\begin{document}

\title{A note on additivity of polygamma functions}

\author[F. Qi]{Feng Qi}
\address[F. Qi]{Research Institute of Mathematical Inequality Theory, Henan Polytechnic University, Jiaozuo City, Henan Province, 454010, China}
\email{\href{mailto: F. Qi <qifeng618@gmail.com>}{qifeng618@gmail.com}, \href{mailto: F. Qi <qifeng618@hotmail.com>}{qifeng618@hotmail.com}, \href{mailto: F. Qi <qifeng618@qq.com>}{qifeng618@qq.com}}
\urladdr{\url{http://qifeng618.spaces.live.com}}

\author[B.-N. Guo]{Bai-Ni Guo}
\address[B.-N. Guo]{School of Mathematics and Informatics,
Henan Polytechnic University, Jiaozuo City, Henan Province, 454010, China}
\email{\href{mailto: B.-N. Guo <bai.ni.guo@gmail.com>}{bai.ni.guo@gmail.com}, \href{mailto: B.-N. Guo <bai.ni.guo@hotmail.com>}{bai.ni.guo@hotmail.com}}
\urladdr{\url{http://guobaini.spaces.live.com}}

\begin{abstract}
In the note, the functions $\abs{\psi^{(i)}(e^x)}$ for $i\in\mathbb{N}$ are proved to be sub-additive on $(\ln\theta_i,\infty)$ and super-additive on $(-\infty,\ln\theta_i)$, where $\theta_i\in(0,1)$ is the unique root of equation $2\abs{\psi^{(i)}(\theta)}=\abs{\psi^{(i)}(\theta^2)}$.
\end{abstract}



\thanks{This paper was typeset using \AmS-\LaTeX}

\maketitle
\section{Introduction}

Recall \cite{alzer-koumandos, alzer-rus, gao-sub-add-psi} that a function $f$ is said to be sub-additive on $I$ if
\begin{equation}\label{sub-dfn-ineq}
f(x+y)\le f(x)+f(y)
\end{equation}
holds for all $x,y\in I$ such that $x+y\in I$. If the inequality \eqref{sub-dfn-ineq} is reversed, then $f$ is called super-additive on $I$.
\par
The sub-additive and super-additive functions play an important role in the
theory of differential equations, in the study of semi-groups, in number
theory, and also in the theory of convex bodies. A lot of literature for the
sub-additive and super-additive functions can be found in \cite{alzer-koumandos,
alzer-rus} and related references therein.
\par
It is well-known that the classical Euler gamma function $\Gamma(x)$ may be defined for $x>0$ by
\begin{equation}\label{egamma}
\Gamma(x)=\int^\infty_0t^{x-1} e^{-t}\td t.
\end{equation}
The logarithmic derivative of $\Gamma(x)$, denoted by $\psi(x)=\frac{\Gamma'(x)}{\Gamma(x)}$, is called the psi or digamma function, and $\psi^{(k)}(x)$ for $k\in\mathbb{N}$ are called the polygamma functions. It is common knowledge that these functions are fundamental and important and that they have much extensive applications in mathematical sciences.
\par
In \cite{Alzer-Ruehr}, the function $\psi(a+x)$ is proved to be sub-multiplicative with respect to $x\in[0,\infty)$ if and only if $a\ge a_0$, where $a_0$ denotes the only positive real number which satisfies $\psi(a_0)=1$.
\par
In~\cite{alzer-rus}, the function $[\Gamma(x)]^\alpha$ was proved to be sub-additive on $(0,\infty)$ if and only if $\frac{\ln2}{\ln\Delta}\le\alpha\le0$, where $\Delta=\min_{x\ge0}\frac{\Gamma(2x)}{\Gamma(x)}$.
\par
In~\cite[Lemma~2.4]{forum-alzer}, the function $\psi(e^x)$ was proved to be strictly concave on $\mathbb{R}$.
\par
In~\cite[Theorem~3.1]{gao-sub-add-psi}, the function $\psi(a+e^x)$ is proved to be sub-additive on $(-\infty,\infty)$ if and only if $a\ge c_0$, where $c_0$ is the only positive zero of $\psi(x)$.
\par
In~\cite[Theorem~1]{Extension-alzer-AMEN.tex}, among other things, it was presented that the function $\psi^{(k)}(e^x)$ for $k\in\mathbb{N}$ is concave $($or convex, respectively$)$ on $\mathbb{R}$ if $k=2n-2$ $($or $k=2n-1$, respectively$)$ for $n\in\mathbb{N}$.
\par
In this short note, we discuss sub-additive and super-additive properties of polygamma functions $\psi^{(i)}(x)$ for $i\in\mathbb{N}$.
\par
Our main result is the following Theorem~\ref{sub-sup-add}.

\begin{thm}\label{sub-sup-add}
The functions $\abs{\psi^{(i)}(e^x)}$ for $i\in\mathbb{N}$ are super-additive on $(-\infty,\ln\theta_i)$ or sub-additive on $(\ln\theta_i,\infty)$, where
$\theta_i\in(0,1)$ is the unique root of equation
\begin{equation}
2\abs{\psi^{(i)}(\theta)}=\abs{\psi^{(i)}(\theta^2)}.
\end{equation}
\end{thm}

\section{Proof of Theorem \ref{sub-sup-add}}

Let
\begin{equation}
f(x,y)=\abs{\psi^{(i)}(x)}+\abs{\psi^{(i)}(y)}-\abs{\psi^{(i)}(xy)}
\end{equation}
for $x>0$ and $y>0$, where $i\in\mathbb{N}$. In order to show
Theorem~\ref{sub-sup-add}, it is sufficient to prove the positivity or
negativity of the function $f(x,y)$. Direct differentiation yields
\begin{equation}
\begin{split}
\frac{\partial f(x,y)}{\partial x} &=y\abs{\psi^{(i+1)}(xy)}-\abs{\psi^{(i+1)}(x)}\\
&=\frac1x\Bigl[xy\abs{\psi^{(i+1)}(xy)}-x\abs{\psi^{(i+1)}(x)}\Bigr].
\end{split}
\end{equation}
In \cite[Lemma~1]{alzeraeq} and~\cite{subadditive-qi.tex, sub-additive-qi-3.tex}, among other things, the functions $x^\alpha\abs{\psi^{(i)}(x)}$ are proved to be strictly increasing on $(0,\infty)$ if and only if $\alpha\ge i+1$ and strictly decreasing if and only if $\alpha\le i$. From this monotonicity, it follows easily that $\frac{\partial f(x,y)}{\partial x}\gtreqless0$ if and only if $y\lesseqgtr1$, which means that the function $f(x,y)$ is strictly increasing for $y<1$ and strictly decreasing for $y>1$ in $x\in(0,\infty)$. Since
$$
\lim_{x\to\infty}f(x,y)=\abs{\psi^{(i)}(y)}>0,
$$
then the function $f(x,y)$ is positive in $x\in(0,\infty)$ for $y>1$.
\par
For $y<1$, by virtue of the increasing monotonicity of $f(x,y)$, it is deduced that
\begin{enumerate}
\item
if $x>1$, then $f(1,y)=\abs{\psi^{(i)}(1)}<f(x,y)<\abs{\psi^{(i)}(y)}$;
\item
if $x<1$, then $f(x,y)<f(1,y)=\abs{\psi^{(i)}(1)}$;
\item
if $y<x<1$, then $f(y,y)<f(x,y)$;
\item
if $x<y<1$, then $f(x,x)<f(x,y)$.
\end{enumerate}
This implies that
\begin{equation}
f(\theta,\theta)=2\abs{\psi^{(i)}(\theta)}-\abs{\psi^{(i)}(\theta^2)}<f(x,y)
\end{equation}
for $y<1$, where $\theta<1$ with $\theta<x$ and $\theta<y$.
Since $f(\theta,\theta)$ is strictly increasing on $(0,1)$ such that $f(1,1)=\abs{\psi^{(i)}(1)}>0$ and $\lim_{\theta\to0^+}f(\theta,\theta)=-\infty$, then the function $f(\theta,\theta)$ has a unique zero $\theta_i\in(0,1)$ such that $f(\theta,\theta)>0$ for $1>\theta>\theta_i$.
\par
In conclusion, the function $f(x,y)$ is positive for $x,y>\theta_i$ or negative for $0<x,y<\theta_i$. The proof of Theorem~\ref{sub-sup-add} is complete.

\end{document}